\documentclass[12pt]{article}%
\usepackage{graphicx}
\usepackage[intlimits]{amsmath}
\usepackage{latexsym}
\usepackage{amsfonts}
\usepackage{amssymb}%
\makeatletter
\newfont{\footsc}{cmcsc10 at 8truept}
\newfont{\footbf}{cmbx10 at 8truept}
\newfont{\footrm}{cmr10 at 10truept}
\newtheorem{theorem}{Theorem}

\newtheorem{proposition}[theorem]{Proposition}

\newenvironment{proof}[1][Proof]{\noindent{\textbf {#1}  }}  {\hfill$\Box$\bigskip}

\begin{document}

\title{Bounds on graph eigenvalues II}
\author{Vladimir Nikiforov\\{\small Department of Mathematical Sciences, University of Memphis, }\\{\small Memphis TN 38152, USA, e-mail: }\textit{vnikifrv@memphis.edu}}
\maketitle

\begin{abstract}
We prove three results about the spectral radius $\mu\left(  G\right)  $ of a
graph $G:$

(a) Let $T_{r}\left(  n\right)  $ be the $r$-partite Tur\'{a}n graph of order
$n.$ If $G$ is a $K_{r+1}$-free graph of order $n,$ then%
\[
\mu\left(  G\right)  <\mu\left(  T_{r}\left(  n\right)  \right)
\]
unless $G=T_{r}\left(  n\right)  .$

(b) For most irregular graphs $G$ of order $n$ and size $m,$
\[
\mu\left(  G\right)  -2m/n>1/\left(  2m+2n\right)  .
\]

(c) Let $0\leq k\leq l.$ If $G$ is a graph of order $n$ with no $K_{2}%
+\overline{K}_{k+1}$ and no $K_{2,l+1},$ then
\[
\mu\left(  G\right)  \leq\min\left\{  \Delta\left(  G\right)  ,\left(
k-l+1+\sqrt{\left(  k-l+1\right)  ^{2}+4l(n-1)}\right)  /2\right\}  .
\]
\bigskip

\textbf{Keywords: }\textit{clique number, spectral radius, Tur\'{a}n graph,
maximum degree, books }

\end{abstract}

\section{Introduction}

Our notation is standard (e.g., see \cite{Bol98}); specifically, we write
$G\left(  n\right)  $ for a graph of order $n$ and $\mu\left(  G\right)  $ for
the maximum eigenvalue of the adjacency matrix of $G.$

Write $T_{r}\left(  n\right)  $ for the $r$-partite Tur\'{a}n graph of order
$n$ and let $G=G\left(  n\right)  .$ In \cite{FLZ07} it is shown that if $G$
is $r$-partite, then $\mu\left(  G\right)  <\mu\left(  T_{r}\left(  n\right)
\right)  $ unless $G=T_{r}\left(  n\right)  $. On the other hand, Wilf
\cite{Wil86} showed that if $G$ is $K_{r+1}$-free, then $\mu\left(  G\right)
\leq\left(  1-1/r\right)  n.$ We strengthen these two results as follows.

\begin{theorem}
\label{th1} If $G=G\left(  n\right)  $ is a $K_{r+1}$-free graph, then
$\mu\left(  G\right)  <\mu\left(  T_{r}\left(  n\right)  \right)  $ unless
$G=T_{r}\left(  n\right)  .$
\end{theorem}

Next, let $G$ be a graph of order $n,$ size $m,$ and maximum degree
$\Delta\left(  G\right)  =\Delta.$ One of the best known facts about
$\mu\left(  G\right)  $ is the inequality $\mu\left(  G\right)  \geq2m/n,$ due
to Collatz and Sinogovitz \cite{CoSi57}. In \cite{Nik06} we gave upper and
lower bounds on $\mu\left(  G\right)  -2m/n$ in terms of degree deviation. In
turn, Cioab\u{a} and Gregory \cite{CiGr07} showed that, if $G$ is irregular
and $n\geq4,$ then $\mu-2m/n>1/\left(  n\Delta+2n\right)  .$ In this note we
give another proof of this bound and improve it for most graphs.

Call a graph \emph{subregular} if $\Delta\left(  G\right)  -\delta\left(
G\right)  =1$ and all but one vertices have the same degree.

\begin{theorem}
\label{th2}If $G$ is an irregular graph of order $n\geq4$ and size $m$, then
\begin{equation}
\mu\left(  G\right)  -2m/n>1/\left(  2m+2n\right)  \label{in1}%
\end{equation}
unless $G$ is subregular. If $G$ is subregular with $\Delta\left(  G\right)
=\Delta$, then
\begin{equation}
\mu\left(  G\right)  -2m/n>1/\left(  n\Delta+2n\right)  . \label{in2}%
\end{equation}

\end{theorem}

Finally, write $B_{k}$ for the graph consisting of $k$ triangles sharing an
edge, and let $0\leq k\leq l\leq\Delta.$ Shi and Song \cite{ShSo07} showed
that if $G=G\left(  n\right)  $ is a connected graph with $\Delta\left(
G\right)  =\Delta,$ with no $B_{k+1}$ and no $K_{2,l+1},$ then
\begin{equation}
\mu\left(  G\right)  \leq\left(  k-l+\sqrt{\left(  k-l\right)  ^{2}%
+4\Delta+4l(n-1)}\right)  /2. \label{in3.0}%
\end{equation}

We extend this result as follows.

\begin{theorem}
\label{th3}Let $0\leq k\leq l.$ If $G=G\left(  n\right)  $ is a graph with
$\Delta\left(  G\right)  =\Delta,$ with no $B_{k+1}$ and no $K_{2,l+1},$ then
\begin{equation}
\mu\left(  G\right)  \leq\min\left\{  \Delta,\left(  k-l+1+\sqrt{\left(
k-l+1\right)  ^{2}+4l(n-1)}\right)  /2\right\}  . \label{in3}%
\end{equation}
If $G$ is connected, equality holds if and only if one of the following
conditions holds:

(i) $\Delta^{2}-\Delta\left(  k-l+1\right)  \leq l(n-1)$ and $G$ is $\Delta$-regular;

(ii) $\Delta^{2}-\Delta\left(  k-l+1\right)  >l(n-1)$ and every two vertices
of $G$ have $k$ common neighbors if they are adjacent, and $l$ common
neighbors otherwise.
\end{theorem}

We note without a proof that (\ref{in3}) implies (\ref{in3.0}).

\section{Proofs}

\begin{proof}
[\textbf{Proof of Theorem \ref{th1}}]Write $k_{r}\left(  G\right)  $ for the
number of $r$-cliques of $G.$ The following result is given in \cite{Nik02}:
if $G$ is $K_{r+1}$-free graph, then
\begin{equation}
\mu^{r}\left(  G\right)  \leq%
{\textstyle\sum\limits_{s=2}^{r}}
\left(  s-1\right)  k_{s}\left(  G\right)  \mu^{r-s}\left(  G\right)  .
\label{polyn}%
\end{equation}
According to a result of Zykov \cite{Zyk49} (see also Erd\H{o}s \cite{Erd62}),
if the clique number of a graph $G$ is $r$, then $k_{s}\left(  G\right)
<k_{s}\left(  T_{r}\left(  n\right)  \right)  $ for every $2\leq s\leq r$,
unless $G=T_{r}\left(  n\right)  .$ Assuming that $G\neq T_{r}\left(
n\right)  $, Zykov's theorem implies that $k_{s}\left(  G\right)
<k_{s}\left(  T_{r}\left(  n\right)  \right)  $ for every $2\leq s\leq r.$
Hence, in view of (\ref{polyn}), we have%
\[
\mu^{r}\left(  G\right)  <%
{\textstyle\sum\limits_{s=2}^{r}}
\left(  s-1\right)  k_{s}\left(  T_{r}\left(  n\right)  \right)  \mu
^{r-s}\left(  G\right)  .
\]
This implies that $\mu\left(  G\right)  <x,$ where $x$ is the largest root of
the equation
\begin{equation}
x^{r}=%
{\textstyle\sum\limits_{s=2}^{r}}
\left(  s-1\right)  k_{s}\left(  T_{r}\left(  n\right)  \right)  x^{r-s}.
\label{Tce}%
\end{equation}
It is known (see, e.g., \cite{CDS80}, p. 74) that (\ref{Tce}) is the
characteristic equation of the Tur\'{a}n graph; so, $\mu\left(  G\right)
<x=\mu\left(  T_{r}\left(  n\right)  \right)  ,$ completing the proof.
\end{proof}

\bigskip

To simplify the proof of Theorem \ref{th2}, we first prove inequality
(\ref{in1}) for two special graphs.

\begin{proposition}
\label{pro1}Inequality (\ref{in1}) holds if $G$ has $n-2$ vertices of degree
$n-1$ and $2$ vertices of degree $n-2$.
\end{proposition}

\begin{proof}
Clearly, $G$ is the complete graph of order $n$ with one edge removed. Using
the theorem of Finck and Grohmann \cite{FiGr65} (see also \cite{CDS80},
Theorem 2.8), we find that
\[
\mu\left(  G\right)  =\frac{n-3+\sqrt{n^{2}+2n-7}}{2}.
\]
Hence, in view of $2m=n^{2}-n-2,$ we obtain,%
\begin{align*}
\mu\left(  G\right)  -\frac{2m}{n}  &  =\frac{\sqrt{n^{2}+2n-7}-\left(
n+1-\frac{4}{n}\right)  }{2}=\frac{4n-8}{n^{2}\left(  \sqrt{n^{2}%
+2n-7}+\left(  n+1-\frac{4}{n}\right)  \right)  }\\
&  >\frac{4n-8}{n^{2}\left(  n+1+\left(  n+1-\frac{4}{n}\right)  \right)
}\geq\frac{2n-4}{n\left(  n^{2}+n-2\right)  }\geq\frac{1}{n^{2}+n-2}=\frac
{1}{2m+2n},
\end{align*}
completing the proof.
\end{proof}

\bigskip

\begin{proposition}
\label{pro2}Inequality (\ref{in1}) holds if $G$ has $n-2$ vertices of degree
$n-2$ and $2$ vertices of degree $n-1$.
\end{proposition}

\begin{proof}
We easily deduce that $n$ is even, say $n=2k,$ and that $G$ is the complement
of a $\left(  k-1\right)  $-matching. Using the theorem of Finck and Grohmann,
we find that
\[
\mu\left(  G\right)  =\frac{n-3+\sqrt{n^{2}-2n+9}}{2}.
\]
Hence, in view of $2m=n^{2}-2n+2,$ we obtain,%
\begin{align*}
\mu\left(  G\right)  -\frac{2m}{n}  &  =\frac{\sqrt{n^{2}-2n+9}-\left(
n-1+\frac{4}{n}\right)  }{2}=\frac{4n-8}{n^{2}\left(  \sqrt{n^{2}%
-2n+9}+\left(  n-1+\frac{4}{n}\right)  \right)  }\\
&  >\frac{4n-8}{n^{2}\left(  n+1+\left(  n-1+\frac{4}{n}\right)  \right)
}=\frac{2n-4}{n\left(  n^{2}+2\right)  }\geq\frac{1}{n^{2}+2}=\frac{1}{2m+2n},
\end{align*}
completing the proof.
\end{proof}

\bigskip

\begin{proof}
[\textbf{Proof of Theorem \ref{th2}}]Set $V=V\left(  G\right)  ,$ $\mu
=\mu\left(  G\right)  ,$ and $\delta=\delta\left(  G\right)  .$ Assume first
that $G$ is not subregular.\bigskip

\emph{Proof of inequality (\ref{in1})}

Our proof is based on Hofmeister's inequality \cite{Hof88}: $\mu^{2}%
\geq\left(  1/n\right)
{\textstyle\sum_{u\in V}}
d^{2}\left(  u\right)  .$\bigskip

\textbf{Case: }$\Delta-\delta\geq2$

In this case we easily see that
\[%
{\textstyle\sum\limits_{u\in V}}
\left(  d\left(  u\right)  -\frac{2m}{n}\right)  ^{2}\geq2>\frac{2m}%
{m+n}+\frac{n}{4\left(  m+n\right)  ^{2}},
\]
and so,
\[
\mu\geq\sqrt{\frac{1}{n}%
{\textstyle\sum\limits_{u\in V}}
d^{2}\left(  u\right)  }=\sqrt{\frac{1}{n}%
{\textstyle\sum\limits_{u\in V}}
\left(  d\left(  u\right)  -\frac{2m}{n}\right)  ^{2}+\frac{4m^{2}}{n^{2}}%
}>\frac{2m}{n}+\frac{1}{2m+2n},
\]
as claimed. Thus, hereafter we shall assume that $\Delta-\delta=1.$\bigskip

\textbf{Case: }$\Delta-\delta=1$

Letting $k$ be the number of vertices of degree $\Delta=\delta+1,$ we have
$2m/n=\delta+k/n,$ and so,
\[
\frac{1}{n}%
{\textstyle\sum\limits_{u\in V}}
\left(  d\left(  u\right)  -\frac{2m}{n}\right)  ^{2}=\frac{n-k}{n}\left(
\frac{k}{n}\right)  ^{2}+\frac{k}{n}\left(  \frac{n-k}{n}\right)  ^{2}%
=\frac{k\left(  n-k\right)  }{n^{2}}.
\]
Hence, if
\begin{equation}
\frac{k\left(  n-k\right)  }{n^{2}}>\frac{2m}{n\left(  m+n\right)  }+\frac
{1}{4\left(  m+n\right)  ^{2}}, \label{in6}%
\end{equation}
then inequality (\ref{in1}) follows as above. Assume for contradiction that
(\ref{in6}) fails.

Suppose first that either $k=2$ or $n-k=2.$ Since (\ref{in6}) fails, we see
that
\begin{align}
2-\frac{4}{n}  &  =\frac{\left(  n-2\right)  2}{n}\leq\frac{k\left(
n-k\right)  }{n}\leq\frac{2m}{m+n}+\frac{n}{4\left(  m+n\right)  ^{2}%
}\nonumber\\
&  =2-\frac{2n}{m+n}+\frac{n}{4\left(  m+n\right)  ^{2}}. \label{in7}%
\end{align}
In view of Propositions \ref{pro1} and \ref{pro2}, we may assume that
$\delta\leq n-3,$ and so%
\[
2m=\delta n+k\leq\delta n+n-2\leq n^{2}-2n-2.
\]
Noting that (\ref{in7}) increases in $m$, we obtain
\[
-\frac{4}{n^{2}}\leq-\frac{4}{n^{2}-2}+\frac{1}{\left(  n^{2}-2\right)  ^{2}%
},
\]
a contradiction for $n\geq4$.

Finally, let $k\geq3$ and $n-k\geq3;$ thus, $n\geq6.$ We have
\[
2m=\delta n+k\leq\delta n+n-3\leq\left(  n-2\right)  n+n-3.
\]
By assumption inequality (\ref{in6}) fails; hence,
\[
3-\frac{9}{n}\leq\frac{\left(  n-k\right)  k}{n}\leq2-\frac{2n}{m+n}+\frac
{n}{4\left(  m+n\right)  ^{2}}\leq2-\frac{4n}{n^{2}+n-3}+\frac{n}{\left(
n^{2}+n-3\right)  ^{2}}.
\]
This inequality is a contradiction for $n\geq6,$ completing the proof of
(\ref{in6}).\bigskip

\emph{Proof of inequality (\ref{in2}) when }$G$\emph{ is subregular }

Since $G$ is subregular, it has either a single vertex of degree $\Delta$ or a
single vertex of degree $\delta.$ Clearly, $\delta\geq1$, and so
$m>n/2.$\bigskip

\textbf{Case: }$G$\textbf{ has a single vertex of degree }$\Delta$

Setting $\Delta=k+1$ and
\[
c=\frac{nk+1}{n}+\frac{1}{n\left(  k+3\right)  }=k+\frac{k+4}{n\left(
k+3\right)  },
\]
in view of $2m=nk+1,$ inequality (\ref{in2}) amounts to $\mu>c\left(
G\right)  .$

Select a vertex $u\in V$ with $d\left(  u\right)  =k+1;$ partition $V$ as
$V=\left\{  u\right\}  \cup V\backslash\left\{  u\right\}  $ and let $B$ be
the quotient matrix of this partition (see, e.g. \cite{GoRo01}, Ch. 9), i.e.,%
\[
B=\left(
\begin{array}
[c]{cc}%
0 & \frac{k+1}{n-1}\\
k+1 & k-\frac{k+1}{n-1}%
\end{array}
\right)  .
\]
Writing $P\left(  x\right)  $ for the characteristic polynomial of $B$ and
observing that $k\leq n-2,$ we have%
\begin{align*}
P\left(  c\right)   &  =\left(  k+\frac{k+4}{n\left(  k+3\right)  }\right)
\left(  k+\frac{k+4}{n\left(  k+3\right)  }-\left(  k-\frac{k+1}{n-1}\right)
\right)  -\frac{\left(  k+1\right)  ^{2}}{n-1}\\
&  =k\frac{k+4}{n\left(  k+3\right)  }+\frac{1}{n^{2}}\left(  \frac{k+4}%
{k+3}\right)  ^{2}+\left(  \frac{k+4}{n\left(  k+3\right)  }\right)
\frac{k+1}{n-1}-\frac{k+1}{n-1}\\
&  =-\frac{3}{n\left(  k+3\right)  }+\frac{1}{n^{2}}\left(  \frac{k+4}%
{k+3}\right)  ^{2}+\frac{\left(  k+1\right)  }{n\left(  n-1\right)  \left(
k+3\right)  }\\
&  =\frac{1}{n^{2}\left(  k+3\right)  }\left(  -3n+2k+6+\frac{1}{k+3}%
+\frac{k+1}{n-1}\right) \\
&  \leq\frac{1}{n^{2}\left(  k+3\right)  }\left(  -3n+2\left(  n-2\right)
+6+\frac{1}{4}+1\right)  <0.
\end{align*}
By interlacing, $P\left(  \mu\right)  \geq0>P\left(  c\right)  ,$ and so
$\mu>c,$ completing the proof of (\ref{in2}) in this case.\bigskip

\textbf{Case: }$G$\textbf{ has a single vertex of degree }$\delta$

Setting $\Delta=k$ and
\[
c=\frac{nk-1}{n}+\frac{1}{n\left(  k+3\right)  }=k-\frac{k+1}{n\left(
k+2\right)  },
\]
in view of $2m=nk-1,$ inequality (\ref{in2}) amounts to $\mu>c.$

Select $u\in V$ with $d\left(  u\right)  =k-1;$ partition $V$ as $V=\left\{
u\right\}  \cup V\backslash\left\{  u\right\}  $ and let $B$ be the quotient
matrix of this partition, i.e.,%
\[
B=\left(
\begin{array}
[c]{cc}%
0 & \frac{k-1}{n-1}\\
k-1 & k-\frac{k-1}{n-1}%
\end{array}
\right)  .
\]
Writing $P\left(  x\right)  $ for the characteristic polynomial of $B$ and
observing that $k\leq n-2,$ we have
\begin{align*}
P\left(  c\right)   &  =\left(  k-\frac{k+1}{n\left(  k+2\right)  }\right)
\left(  k-\frac{k+1}{n\left(  k+2\right)  }-\left(  k-\frac{k-1}{n-1}\right)
\right)  -\frac{\left(  k-1\right)  ^{2}}{n-1}\\
&  =-\frac{k\left(  k+1\right)  }{n\left(  k+2\right)  }+\frac{1}{n^{2}%
}\left(  \frac{k+1}{k+2}\right)  ^{2}+\frac{k-1}{n\left(  n-1\right)  \left(
k+2\right)  }+\frac{k-1}{n}\\
&  =-\frac{2}{n\left(  k+2\right)  }+\frac{1}{n^{2}}\left(  \frac{k+1}%
{k+2}\right)  ^{2}+\frac{k-1}{n\left(  n-1\right)  \left(  k+2\right)  }\\
&  =\frac{1}{n^{2}\left(  k+2\right)  }\left(  -2n+2k+1+\frac{1}{k+2}%
+\frac{k-1}{\left(  n-1\right)  }\right) \\
&  <\frac{1}{n^{2}\left(  k+2\right)  }\left(  -2n+2\left(  n-2\right)
+1+\frac{1}{1+2}+1\right)  <0.
\end{align*}
By interlacing, $P\left(  \mu\right)  \geq0>P\left(  c\right)  ,$ completing
the proof of (\ref{in2}).\bigskip
\end{proof}

\begin{proof}
[\textbf{Proof of Theorem \ref{th3}}]Set $V=V\left(  G\right)  $ and $\mu
=\mu\left(  G\right)  ;$ given $u\in V,$ write $\Gamma\left(  u\right)  $ for
the set of neighbors of $u.$ Select $u\in V;$ let $A=\Gamma\left(  u\right)
$, $B=V\backslash\left(  \Gamma\left(  u\right)  \cup\left\{  u\right\}
\right)  ,$ and $e\left(  A,B\right)  $ be the number of $A-B$ edges. Since
$G$ contains no $B_{k+1}$ and no $K_{2,l+1},$ we see that%
\begin{equation}%
{\textstyle\sum\limits_{v\in A}}
\left(  d\left(  v\right)  -k-1\right)  \leq%
{\textstyle\sum\limits_{v\in A}}
\left\vert \Gamma\left(  v\right)  \cap B\right\vert =e\left(  A,B\right)  =%
{\textstyle\sum\limits_{v\in A}}
\left\vert \Gamma\left(  v\right)  \cap A\right\vert \leq\left(  n-d\left(
u\right)  -1\right)  l. \label{in5}%
\end{equation}

Letting $A$ be the adjacency matrix of $G,$ note that the $u$th row sum of the
matrix
\[
C=A^{2}-\left(  k+1-l\right)  A-\left(  n-1\right)  lI_{n}%
\]
is equal to
\[%
{\textstyle\sum\limits_{v\in A}}
\left(  d\left(  v\right)  -k-1\right)  -\left(  n-1-d\left(  u\right)
\right)  l;
\]
consequently, all row sums $C$ are nonpositive. Letting $\mathbf{x}=\left(
x_{1},\ldots,x_{n}\right)  $ be an eigenvector of $A$ to $\mu,$ we see that
the value
\[
\lambda=\mu^{2}-\left(  k+1-l\right)  \mu-\left(  n-1\right)  l
\]
is an eigenvalue of $C$ with eigenvector $\mathbf{x}$. Therefore, $\lambda
\leq0,$ and so,%
\[
\mu\leq\left(  k-l+1+\sqrt{\left(  k-l+1\right)  ^{2}+4l(n-1)}\right)  /2,
\]
completing the proof of inequality (\ref{in3}).

Let equality hold in (\ref{in3}) and $G$ be connected; thus, the eigenvector
$\mathbf{x}=\left(  x_{1},\ldots,x_{n}\right)  $ to $\mu$ is positive. We
shall prove the necessity of conditions \emph{(i)} and \emph{(ii)}. If
\[
\mu=\Delta\leq\left(  k-l+1+\sqrt{\left(  k-l+1\right)  ^{2}+4l(n-1)}\right)
/2,
\]
then $\Delta^{2}-\Delta\left(  k-l+1\right)  \leq l(n-1)$ and $G$ is $\Delta$-regular.

On the other hand, if
\[
\mu=\left(  k-l+1+\sqrt{\left(  k-l+1\right)  ^{2}+4l(n-1)}\right)
/2<\Delta,
\]
then $\Delta^{2}-\Delta\left(  k-l+1\right)  >l(n-1)$ and $\lambda=0.$ Scaling
$\mathbf{x}$ so that $x_{1}+\cdots+x_{n}=1,$ we see that $\lambda$ is a convex
combination of the row sums of $C$ which are nonpositive; thus, all row sums
of $C$ are $0.$ Since equality holds in (\ref{in5}) for every $u\in\left[
n\right]  ,$ every two vertices have exactly $k$ common neighbors if they are
adjacent, and exactly $l$ common neighbors otherwise. This completes the proof.
\end{proof}

\section{\textbf{Concluding remarks}}

Finding tight bounds on the spectral radius of subregular graphs is a
challenging problem. Specifically, we cannot determine for which subregular
graphs $G$ one has
\[
\mu\left(  G\right)  >\frac{2m}{n}+\frac{1}{2m+2n}.
\]

Note that strongly regular graphs satisfy condition \emph{(ii)} for equality
in (\ref{in3}), but irregular graphs can satisfy this condition as well, e.g.,
the star $K_{1,n-1}$ and the friendship graph.

Finally, setting $l=\Delta$ or $k=0,$ Theorem \ref{th3} implies assertions
that strengthen Corollaries 1 and 2 of \cite{ShSo07}.


\begin{thebibliography}{99}                                                                                               %


\bibitem {Bol98}B. Bollob\'{a}s, \emph{Modern Graph Theory}\textit{,} Graduate
Texts in Mathematics, \textbf{184}, Springer-Verlag, New York (1998), xiv+394 pp.

\bibitem {CDS80}D. Cvetkovi\'{c}, M. Doob, and H. Sachs, \emph{Spectra of
Graphs,} VEB Deutscher Verlag der Wissenschaften, Berlin, 1980, 368 pp.

\bibitem {CiGr07}S. M. Cioab\u{a}, D. Gregory, Large matchings from
eigenvalues, to appear in \emph{Linear Algebra Appl. }

\bibitem {CoSi57}L. Collatz, U. Sinogowitz, Spektren endlicher Grafen,
\emph{Abh. Math. Sem. Univ. Hamburg} \textbf{21} (1957), 63-77.

\bibitem {Erd62}P. Erd\H{o}s, On the number of complete subgraphs contained in
certain graphs, \emph{Publ. Math. Inst. Hung. Acad. Sci.} \textbf{VII}, Ser.
A3 (1962), 459-464p.

\bibitem {FiGr65}H.J.Finck, G. Grohmann, Vollst\"{a}ndiges Produkt,
chromatische Zahl und charakteristisches Polynom regul\"{a}rer Graphen. I.
(German) \emph{Wiss. Z. Techn. Hochsch. Ilmenau} 11 (1965) 1--3.

\bibitem {FLZ07}L. Feng, Q. Li, X.-D. Zhang, Spectral radii of graphs with
given chromatic number, \emph{Appl. Math. Letters}, \textbf{20} (2007), 158-162.

\bibitem {GoRo01}C. Godsil, G. Royle, \emph{Algebraic graph theory}, Graduate
Texts in Mathematics, \textbf{207}, Springer-Verlag, New York (2001), xx+439 pp.

\bibitem {Hof88}M. Hofmeister, Spectral radius and degree sequence,
\emph{Math. Nachr.} \textbf{139 }(1988), 37-44.

\bibitem {Nik02}V. Nikiforov, Some inequalities for the largest eigenvalue of
a graph. \emph{Combin. Probab. Comput.} \textbf{11} (2002), 179--189.

\bibitem {Nik06}V. Nikiforov, Eigenvalues and degree deviation in graphs,
\emph{Linear Algebra Appl. }\textbf{414} (2006), 347-360.

\bibitem {ShSo07}L. Shi, Z. Song, Upper bounds on the spectral radius of
book-free and/or $K_{2,l+1}$ graphs, \emph{Linear Algebra Appl. }\textbf{420}
(2007), 526--529.

\bibitem {Wil86}H. Wilf, Spectral bounds for the clique and independence
numbers of graphs, J. Combin. Theory Ser. B 40 (1986), 113-117.

\bibitem {Zyk49}A. A. Zykov, On some properties of linear complexes (in
Russian), \emph{Mat. Sbornik N.S.} \textbf{24}(66), (1949), 163--188.
\end{thebibliography}
\end{document}